# OBJECT ORIENTED DATA ANALYSIS: SETS OF TREES

By Haonan Wang[1] and J. S. Marron[1,2]

*Colorado State University and University of North Carolina*

Object oriented data analysis is the statistical analysis of populations of complex objects. In the special case of functional data analysis, these data objects are curves, where standard Euclidean approaches, such as principal component analysis, have been very successful. Recent developments in medical image analysis motivate the statistical analysis of populations of more complex data objects which are elements of mildly non-Euclidean spaces, such as Lie groups and symmetric spaces, or of strongly non-Euclidean spaces, such as spaces of tree-structured data objects. These new contexts for object oriented data analysis create several potentially large new interfaces between mathematics and statistics. This point is illustrated through the careful development of a novel mathematical framework for statistical analysis of populations of tree-structured objects.

**1. Introduction.** Object oriented data analysis (OODA) is the statistical analysis of data sets of complex objects. The area is understood through consideration of the *atom of the statistical analysis*. In a first course in statistics, the atoms are numbers. Atoms are vectors in multivariate analysis. An interesting special case of OODA is functional data analysis, where atoms are curves; see Ramsay and Silverman [16, 17] for excellent overviews, as well as many interesting analyses, novel methodologies and detailed discussion. More general atoms have also been considered. Locantore et al. [13] studied the case of images as atoms, and Pizer, Thall and Chen [15] and Yushkevich et al. [23] took the atoms to be shape objects in two- and three-dimensional space.

An important major goal of OODA is understanding *population structure* of a data set. The usual first step is to find a *centerpoint*, for example, a

Received January 2005; revised December 2006.

[1]Supported in part by NSF Grant DMS-99-71649 and by NIH Grant NCI P01 CA047982.

[2]Supported by NSF Grant DMS-03-08331.

mean or median, of the data set. The second step is to analyze the *variation about the center*. Principal component analysis (PCA) has been a workhorse method for this, especially when combined with new visualizations as done in functional data analysis. An important reason for this success to date is that the data naturally lie in Euclidean spaces, where standard vector space analyses have proven to be both insightful and effective.

Medical image analysis is motivating some interesting new developments in OODA. These new developments are not in traditional imaging areas, such as the denoising, segmentation and/or registration of a single image, but instead are about the analysis of *populations of images*. Again common goals include finding centerpoints and variation about the center, but also discrimination, that is, classification, is important. A serious challenge to this development is that the data often naturally lie in non-Euclidean spaces. A range of such cases has arisen, from mildly non-Euclidean spaces, such as Lie groups and Riemannian symmetric spaces (see Wikipedia [21, 22] for a good introduction to these concepts), to strongly non-Euclidean spaces, such as populations of tree- or graph-structured data objects. Because such non-Euclidean data spaces are generally unfamiliar to statisticians, there is opportunity for the development of several types of new interfaces between statistics and mathematics. One purpose of this paper is to highlight some of these. The newness of this nonstandard mathematical statistics that is currently under development (and much of which is yet to be developed) is underscored by a particularly deep look at an example of tree-structured data objects.

Lie groups and symmetric spaces are the natural domains for the data objects which arise in the medial representation of body parts, as discussed in Section 1.1. Human organs are represented using vectors of parameters, which have both real valued and angular components. Thus each data object is usefully viewed as a point in a Lie group, or a symmetric space, that is, a curved manifold space. Such representations are often only mildly non-Euclidean, because these curved spaces can frequently be approximated to some degree by tangent spaces, where Euclidean methods of analysis can be used. However the most natural and convincing analysis of the data is done "along the manifold," as discussed in Section 1.1. Because there already exists a substantial medical imaging literature on this, only an overview is given here.

Data objects which are trees or graphs are seen in Section 1.2 to be important in medical image analysis for several reasons. These data types present an even greater challenge, because the data space is strongly non-Euclidean. Fundamental tools of standard vector space statistical analysis, such as linear subspace, projection, analysis of variance and even linear combination are no longer available. Preliminary ad hoc attempts made by the authors at this type of OODA ended up collapsing in a mass of contradictions, because they were based on trying to apply Euclidean notions in this very non-Euclidean domain. This motivated the development of a really new type of mathematical statistics: a rigorous definition–theorem–proof framework for the analysis of such data, which was the dissertation of Wang [19]. In Section 2 it is seen how these tools provide an analysis of a real data set. Section 3 gives an overview of the mathematical structure that underpins the analysis.

Note that statistics and mathematics (of some nonstandard types) meet each other in several ways in OODA. For the Lie group—symmetric space data, mathematics provides a nonstandard framework for conceptualizing the data. For data as trees, an axiomatic system is used as a device to overcome our poor intuition for data analysis in this very non-Euclidean space. Both of these marriages of mathematics and statistics go in quite different directions from that of much of mathematical statistics: the validation and comparison of existing statistical methods through asymptotic analysis as the sample size tends to infinity. Note that this latter type of analysis has so far been completely unexplored for these new types of OODA, and it also should lead to the development of many more interesting connections between mathematics and statistics.

1.1. *OODA on Lie groups—symmetric spaces.* Shape is an interesting and useful characteristic of objects (usually in three dimensions) in medical image analysis. Shape is usually represented as a vector of measurements, so that a data set of shapes can be analyzed as a set of vectors. There are a number of ways to represent shapes of objects. The best known in the statistical literature is landmark based approaches; see Dryden and Mardia [4] for good overview of this area. While they have been a workhorse for solving a wide variety of practical problems, landmark approaches tend to have limited utility for population studies in medical imaging, because a sufficient number of well-defined, replicable landmarks are frequently impossible to define.

Other approaches to shape representation are discussed in Section 1.1 of Wang and Marron [20], but are not discussed here to save space.

A class of convenient and powerful shape representations is *m-reps* (a shortening of "medial representation"), which are based on *medial* ideas; see Pizer, Thall and Chen [15] and Yushkevich et al. [23] for detailed introduction and discussion. The main idea is to find the "central skeletons" of objects, and then to represent the whole object in terms of "spokes" from the center to the boundary. The central structure and set of spokes to the boundary are discretized and approximated by a finite set of m-reps. The m-rep parameters (location, radius and angles) are the features and are concatenated into a feature vector to provide a numerical summary of the shape. Each data object is thus represented as the direct product (thus a large vector) of these parameters over the collection of m-reps. A major motivation for using m-reps over other types of representation is that they provide a more direct solution to the *correspondence problem*, which is to match parts of one object with corresponding parts of other members of the population.

A simple example of the use of m-reps in OODA is shown in Figure 1, which uses the specific representation of Yushkevich et al. [23], which studied a set of human corpora callosa gathered from two-dimensional magnetic resonance images. The corpus callosum is the small window between the left and right halves of the brain. The left hand panel of Figure 1 shows a single m-rep decomposition of one corpus callosum. Each large central dot shows the center of an m-rep (five of which are used to represent this object). The m-reps are a discretization of the medial axis, shown in gray in the center of the object. The boundary of the object is determined by the spokes, which are the shorter line segments emanating from each m-rep. These spokes are paired, and are determined by their (common) angle from the medial axis, and their length. All of these parameters are summarized into a feature vector which is used to represent each object.

The right-hand panel of Figure 1 shows a simple OODA of a population of 72 corpora callosa. This is done here by simple principal component analysis of the set of feature vectors. The central shape is the mean of the population. The sequence of shapes gives insight into population variation, by showing the first principal component (thus the mode of maximal variation). In particular, each shape shows a location along the eigenvector of the first PC. This shows that the dominant mode of variation in this population is in the direction of more overall bending in one direction, versus less overall bending in the opposite direction.

For a much deeper and more complicated example of using m-reps for shape analysis, see Figure 2 of Wang and Marron [20]. A simple approach to OODA for m-rep objects is to simply use Euclidean PCA on the vectors of parameters. However, there is substantial room for improvement, because some parameters are angles, while others are radii (thus positive in sign), and still others are position coordinates. One issue that comes up is that units are not commensurate, so some vector entries could be orders of magnitude different from the others, which will drastically affect PCA. An approach to this problem is to replace the eigen-analysis of the covariance matrix (conventional PCA) with the an eigen-analysis of the correlation matrix (a well-known scale free analog of PCA). But this still does not address the central challenge of statistical analysis of angular data. For example, what is the average of a set of angles where some are just above $0°$, and the rest are just below $360°$? The sensible answer is something very close to $0°$, but simple averaging of the numbers involved can give a diametrically opposite answer closer to $180°$. There is a substantial literature on the statistical analysis of angular data, also called directional data, that is, data on the circle or sphere. See Fisher [5], Fisher, Lewis and Embleton [6], Mardia [10] and Mardia and Jupp [11] for good introduction to this area. A fundamental concept of this area is that the most convenient mathematical representation of angles is as points on the unit circle, and for angles in 3d, as points on the unit sphere.

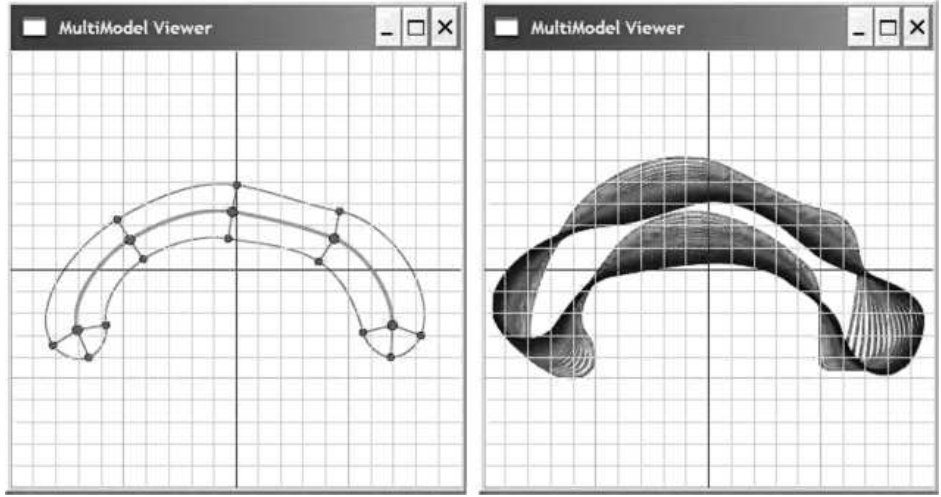

FIG. 1. *Corpus callosum data. Left-hand panel shows all components of the medial representation of the corpus callosum of one person. Right-hand panel shows the boundaries of objects lying along the first PCA eigenvectors, showing the largest component of variation in the population.*

For these same reasons it is natural to represent the vectors of m-rep parameters as direct products of points on the circle and/or sphere for the angles, as positive reals for the radii, and as real numbers for the locations. As noted in Fletcher et al. [7], the natural framework for understanding this type of data object is Lie groups and/or symmetric spaces. Fletcher et al. [7] go on to develop an approach to OODA for such data. The Fréchet approach gives a natural definition of the sample center, and principal geodesic analysis (PGA) quantifies population variation.

The Fréchet mean has been a popular concept in robustness, since it provides useful generalizations of the sample mean. It also provides an effective starting point for non-Euclidean OODA. The main idea is that one way to characterize the sample mean is as the minimizer of the sum of squared distances to each data point. Thus the Fréchet mean can be defined in quite abstract data spaces, as long as a suitable metric can be found. For Lie group–symmetric space data, the natural distance is along geodesics, that is, along the manifold, and this Fréchet mean is called the geodesic mean.

Fletcher's Lie group–symmetric space variation of PCA is PGA. The key to this approach is to characterize PCA as finding lines which maximally approximate the data. On curved manifolds, the analog of lines are geodesics, so PGA searches for geodesics which maximally approximate the data. See Fletcher et al. [7] for detailed discussion and insightful examples.

For other examples of OODA, and discussion of the relationships between OODA and a variety of other areas, see Wang and Marron [20].

1.2. *OODA on tree spaces.* A type of data space which is much farther from Euclidean in nature is the set of trees. A simple motivating example of trees as data is the case of multifigural objects, of the type shown in Figure 2 of Wang and Marron [20]. In that example, all three figures are present in every data object. But if some figures are missing, then the usual vector of m-rep parameters has missing values. Thus the natural data structure is trees, with *nodes* representing the figures. For each figure, the corresponding m-rep parameters appear as *attributes* of that node. A more complex and challenging example is the case of blood vessel trees, discussed in Section 2.

In most of the rest of this paper, the focus is on the very challenging problem of OODA for data sets of tree-structured objects. Tree-structured data objects are mathematically represented as simple graphs (a collection of *nodes* and *edges*, each of which connects some pair of nodes). Simple graphs have a unique path (a set of edges) between every pair of nodes (vertices). A tree is a simple graph, where one node is designated as the *root node*, and all other nodes are *children* of a *parent* node that is closer to the root, where parents and children are connected by edges. In many applications, a tree-structured representation of each data object is very natural, including medical image analysis, phylogenetic studies, clustering analysis and some forms of classification (i.e., discrimination). Limited discussion with references of these areas is given in Section 1.2.1. Our driving example, based on a data set of tree-structured blood vessel trees, is discussed in Section 2.

For a data set of tree-structured data objects, it is unclear how to develop notions such as *centerpoint* and *variation about the center*. Our initial ad hoc attempts at this were confounded by the fact that our usual intuitive ideas led to contradictions. As noted above, we believe this is because our intuition is based on Euclidean ideas, such as linear subspaces, projections and so on, while the space of trees is very "non-Euclidean" in nature, in the sense that natural definitions of the fundamental linear operators of addition and scalar multiplication operations do not seem to be available. Some additional mathematical basis for the claim of "non-Euclideanness of tree space," in the context of phylogenetic trees, can be found in Billera, Holmes and Vogtmann [2]. This failure of our intuition to give the needed insights, has motivated our development of the careful axiomatic mathematical theory for the statistical analysis of data sets of trees given in Section 3. Our approach essentially circumvents the need to define the linear operations that are the foundations of Euclidean space.

The development essentially starts from a Fréchet approach, which is based on a metric. In general, we believe that different data types, such as those listed in Section 1.2.1, will require careful individual choice of a metric. In Sections 3.2 and 3.3, we define a new metric which makes sense for our driving problem of a data set of blood vessel trees.

Once a metric has been chosen, the Fréchet mean of a data set is the point which minimizes the sum of the squared distances to the data points. A simple example is the conventional sample mean in Euclidean space (just the mean vector), which is the Fréchet mean with respect to Euclidean distance. In Section 3.4 this idea is the starting point of our development of the notion of centerpoint for a sample of trees.

After an appropriate centerpoint is defined, it is of interest to quantify the variability of the sample about this center. Here, an analog of PCA, based on the notion of a *treeline* which plays the role of "one-dimensional subspace," is developed for tree space (see Section 3.5). A key theoretical contribution is a fundamental theory of *variation decomposition in tree space*, a tree version of the Pythagorean theorem (see Section 3.5), which allows ANOVA style decomposition of sums of squares.

The driving problem in this paper is the analysis of a sample of blood vessel trees in the human brain; see Bullitt and Aylward [3]. We believe that similar methods could be used for related medical imaging problems, such as the study of samples of pulmonary airway systems, as studied in Tschirren et al. [18]. The blood vessel systems considered here are conveniently represented as trees. In our construction of these trees, each node represents a blood vessel, and the edges only illustrate the connectedness property between two blood vessels. For these blood vessel trees, both topological structure (i.e., connectivity properties) and geometric properties, such as the locations and orientations of the blood vessels, are very important. These geometric properties are summarized as the attributes of each node.

Focusing on our driving example of blood vessel trees and their corresponding attributes, we develop a new metric $\delta$ on tree space; see Section 3.3. Margush [14] gives a deeper discussion of metrics on trees. This metric $\delta$ consists of two parts: the integer part $d_I$, which captures the topological aspects of the tree structure (see Section 3.2 for more detail), and the fractional part $f_\delta$, which captures characteristics of the nodal attributes (see Section 3.3).

The metric $\delta$ provides a foundation for defining the notion of centerpoint. A new centerpoint, the median-mean tree, is introduced (see Section 3.4). It has properties similar to the median with respect to the integer part metric (see Section 3.2) and similar to the mean with respect to the fractional part metric (see Section 3.3).

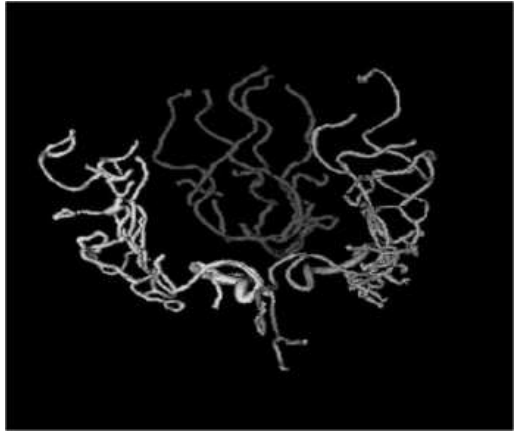

FIG. 2. *The three component blood vessel trees from one person. This detailed graphical illustration uses all attributes from the raw data.*

In Section 3, methods are developed for the OODA of samples of trees. An interesting question for future research is how our sample centerpoint and measures of variation about the center correspond to theoretical notions of these quantities and an underlying probabilistic model for the population. For a promising approach to this problem, see Larget, Simon and Kadane [12].

1.2.1. *Additional applications of OODA for trees.* Our driving application of OODA for tree-structured data objects, to analyze a data set of blood vessel trees, is discussed in Section 2. A number of additional important potential applications, which have not been tried yet, are discussed in Wang and Marron [20].

**2. Tree OODA of a blood vessel data set.** In this section, advanced statistical analysis, including centerpoint and variation about the center, of a data set of tree-structured objects is motivated and demonstrated in the context of human brain blood vessel trees.

An example of arterial brain blood vessels from one person, provided by E. Bullitt, is shown in Figure 2. Because of the branching nature of blood vessel systems, a tree-structured data representation is very natural. See Bullitt and Aylward [3] for detailed discussion of the data, and the method that was used to extract trees of blood vessel systems from magnetic resonance images. The blood vessel systems considered here have three important components: the left carotid, the right carotid and the vertebrobasilar systems, shown in different gray levels in Figure 2. Each component consists of one root vessel and many offspring branches (vessels). Each branch is represented as a node in the tree structure. The attributes for each node include both information about that vessel, and also tree connectivity information. The individual information about that branch is coded as a sequence of vessel medial points (essentially a discretization of the medial axis of the blood vessel), where each point has a 3d location and a radius (of the vessel at the point). The connectivity information for each branch (node) records an index of its parent, and also the location of attachment to the parent. All of these attributes are used in the visual rendering shown in Figure 2.

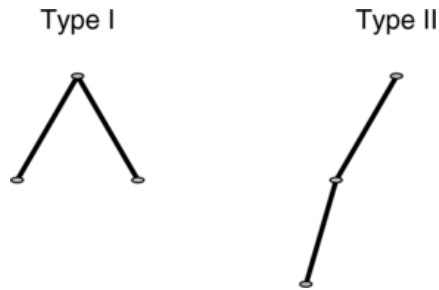

Fig. 3. *Two types of three-node tree structures, that are present in our sample, of simplified blood vessel trees, where seven are of Type* I*, and four are of Type* II.

The full data set analyzed here has 11 trees from three people. These are the left carotid, right carotid and vertebrobasilar systems from each person, plus two smaller, unattached, components from one of the three people.

For simplicity of analysis, in this paper we will work with only a much smaller set of attributes, based on a simple linear approximation of each branch. In particular, the attributes of the root node are the 3d locations of the starting and ending medial points. The attributes of the other branches include the index of the parent, together with a connectivity parameter indicating location of the starting point on the linear approximation of the parent, as

$$p = \frac{\text{Distance of starting point to point of attachment on the parent}}{\text{Distance of starting point to ending point on the parent}},$$

and the 3d locations of the ending point. An additional simplification is that radial information is ignored.

For computational speed, only a subtree (up to three levels and three nodes) of each element among those 11 trees is considered. There are only two different tree structures in this data set, which are called Type I and Type II, shown in Figure 3. Among these 11 blood vessel trees, seven trees have Type I structure and four trees have Type II structure.

Each panel of Figure 4 shows the individual component trees for one person. The three-dimensional aspect of these plots is most clearly visible in rotating views, which are Internet available from the links "first person," "second person" and "third person" on the web site of Wang, www.stat.colostate.edu/~wanghn/tree.htm. These components are shown as thin line trees, which represent each raw data point. Trees are shown using the simplified rendering, based on only the linear approximation attributes, as described above. The root node of each tree is indicated with a solid line type, while the children are dashed. We will first treat each person's component trees as a separate subsample. Each panel of Figure 4 also includes the new notion of *centerpoint* (for that subsample), shown using a thicker line type. This is the *median-mean tree*, as developed in Section 3.4. This tree is central in terms of structure, size and location, in senses which will be defined there.

These trees are combined into a single, larger sample in Figure 6 of Wang and Marron [20]. It turns out that the median-mean tree is surprisingly small, especially in comparison to the median-mean trees for individual people, shown in Figure 4. This will be explained through a careful analysis of the variation about the median-mean tree.

Another important contribution of this paper is the development of an approach to analyzing the variation within a sample of trees. In conventional multivariate analysis, a simple first-order linear approach to this problem is principal component analysis. We develop an analog for samples of trees in Section 3.5. Our first approach is illustrated in Figure 5, with an analysis of the dominant mode of tree structure variation for the full blood vessel tree sample shown in Figure 6 of Wang and Marron [20].

The generalization of PCA to samples of tree-structured objects could be approached in many ways, because PCA can be thought of in a number of different ways. After considering many approaches, we found a suggestion by J. O. Ramsay to be the most natural. The fundamental idea is to view PCA as a sequence of one-dimensional representations of the data. Hence, our tree version PCA is based on notions of *one-dimensional representation* of the data set. These notions are carefully developed and precisely defined in Section 3.5. The foundation of this approach is the concept of *treeline*, which plays the role of line (a one-dimensional subspace in Euclidean space) in tree space. Two different types of treelines are developed in Section 3.5. The *structure treeline* which quantifies sample variation in tree structure, is formally defined in Definition 3.1 and is illustrated here in Figures 5 and 6. The *attribute treeline* describes variation within a fixed type of tree structure, is defined in Definition 3.2, and is illustrated here in Figure 7.

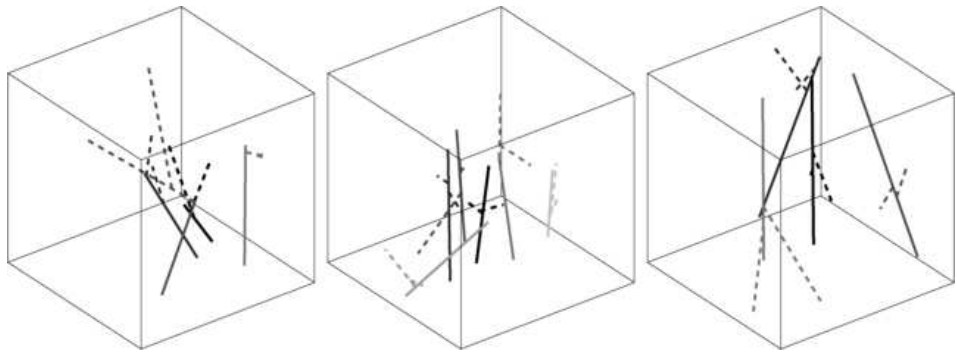

Fig. 4. *Simplified blood vessel trees (thin lines) for each person individually, with the individual median-mean trees (thicker black line). Root nodes use solid line types and children are dashed.*

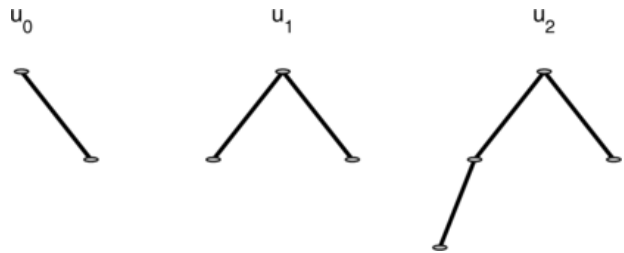

Fig. 5. *The principal structure treeline for the full simplified blood vessel data, without nodal attributes. The figure shows that the dominant sample variation in structure is towards the addition of left-hand children nodes.*

The structure treeline which best represents the data set (this will be formally defined in Section 3.5, but for now think in analogy to PCA) is called the *principal structure treeline*. The principal structure treeline for the full simplified blood vessel data is shown in Figure 5 (structure only, without attributes) and Figure 6 (with attributes). In Figure 5 this treeline starts with the tree $u_0$, which has two nodes. The other trees in this direction are $u_1$ and $u_2$, which consecutively add one left child. Generally structure treelines follow the pattern of successively adding single child nodes. This principal structure treeline is chosen, among all treelines that pass through the median-mean tree, to explain as much of the structure in the data as possible (in a sense defined formally in Section 3.5). Hence, this highlights structure variation in this sample by showing that the dominant component of topological structure variation in the data set is toward branching in the direction of addition of left-hand children nodes. Next, we also study how the attributes change as we move along this principal structure treeline, in Figure 6. The three panels show the simplified tree rendering of the trees whose structure is illustrated in Figure 5, with the first treeline member $u_0$ shown in the left box, the three node tree $u_1$, which is the median-mean tree, in the center box, and $u_2$ with four nodes in the right hand box.

In addition to the principal structure representation, another useful view of the data comes from the principal attribute directions (developed in Section 3.5). *Principal attribute treelines* have a fixed tree structure, and highlight important sample variation within the given tree structure. Since the tree structure is fixed, this treeline is quite similar to the conventional first principal component, within that structure. Here we illustrate this idea, showing the principal attribute treeline which passes through (and thus has the same structure as) the median-mean tree, shown in Figure 7. There are six subplots in this figure. The subplots depict a succession of locations on the attribute treeline, which highlights the sample variation in this treeline direction. These are snapshots which are extracted from a movie version that provides clear visual interpretation of this treeline, and is Internet available from the link "median-mean tree" from Wang's website, **www.stat.colostate.edu/wanghn/tree.htm**. A similar movie, showing a different principal attribute treeline can be found at the link "support tree" (this concept is explained in Section 3.1) on the website.

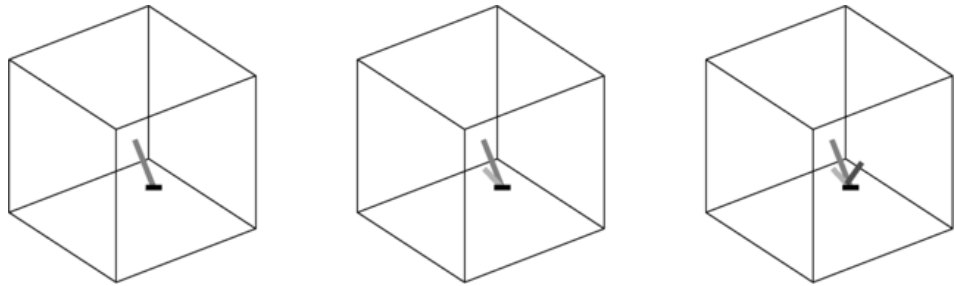

Fig. 6.    *The principal structure treeline with nodal attributes. This shows more about the sample variation than is available from mere structure information.*

In general, Figure 7 shows marked change in the length and orientation of the main root (solid black line). It starts (upper left) as a long nearly vertical segment, which becomes shorter and moves towards horizontal (upper right). This trend continues in the lower left box, where the root is very short indeed and is horizontal. In the next plots (lower row) the root begins to grow, this time in the opposite direction. In particular, the root node flips over, with the top and bottom ends trading places. While these trends are visible here, the impression is much clearer in the movie version. The branches also change in a way that shows smaller scale variation in the data. This was a surprising feature of the sample. Careful investigation showed that the given data sets did not all correctly follow the protocol of choosing the coordinate system according to the direction of blood flow. Some of them have the same direction, while some of them have the reverse direction. A way of highlighting the two different data types is via the *projections* (the direct analogs of the principal component coefficients in PCA) of the 11 trees on this attribute treeline, as shown in Figure 10 of Wang and Marron [20]. This shows that there are two distinct clusters with a clear gap in the middle, validating the above discovery. The gap also shows that no trees correspond to the fourth frame in Figure 7, with a very short root, which can also be seen in the raw data in Figure 6 of Wang and Marron [20]. This shows that the surprisingly short root node, for the median-mean tree, resulted from its being central to the sample formed by these two rather different subgroups that were formed by different orientations of the blood flow in the data set. Note that the clusters dominate the total variation, perhaps obscuring population features of more biological interest.

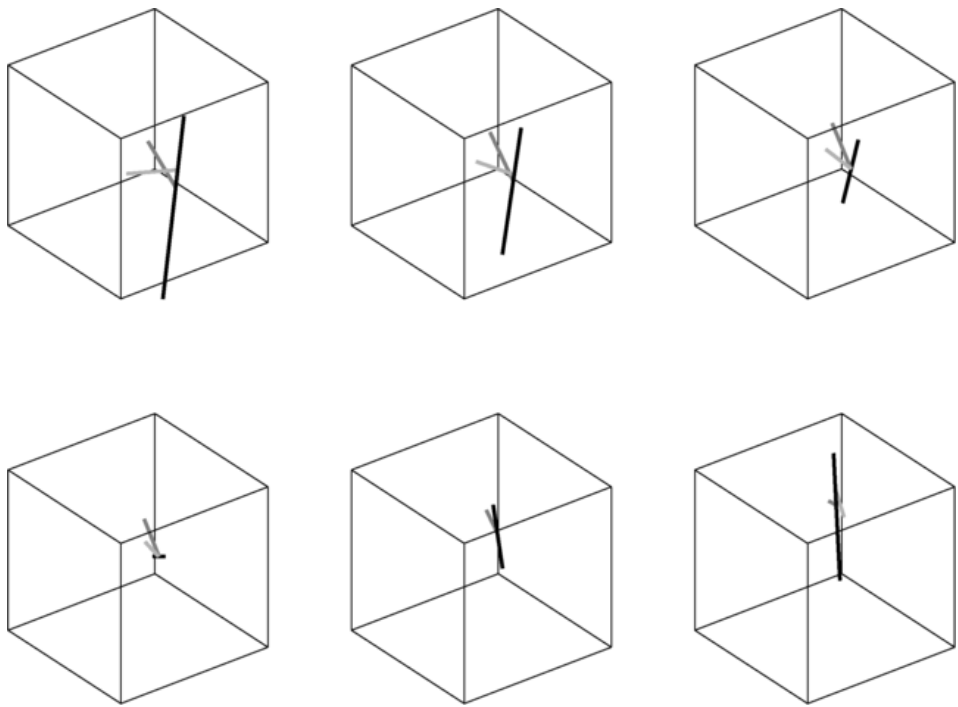

Fig. 7. *The principal attribute treeline passing through the median-mean tree. These are snapshots from a movie which highlights this mode of variation in the sample. The thick black root node flips over.*

**3. Development of the tree OODA methodology.** In this section a rigorous mathematical foundation is developed for the OODA of a data set of trees. We will use $\mathcal{S} = \{t_1, t_2, \ldots, t_n\}$ to denote the data set of size $n$. Careful mathematics are needed because the non-Euclidean nature of tree space means that many classical notions do not carry over as expected. For simplicity, only the case of *binary* trees with finite level is explicitly studied. A binary tree is a tree such that every node has at most two children (left child and right child). If a node has only one child, it should be designated as one of left and right. In our blood vessel application, we consistently label each single child as left. The set of all binary trees, the *binary tree space*, is denoted by $\mathcal{T}$.

3.1. *Notation and preliminaries.* This section introduces a labeling system for the nodes of each tree in the sample, that is, each $t \in \mathcal{S}$. Each tree has a designated node called the *root*. An important indicator of node location in the tree is the *level of the node*, which is the length (number of edges) of the path to the root. In addition, it is convenient to uniquely label each node of a binary tree by a natural number, called the *level-order index*. The level-order index of the node $\omega$ is denoted by $\mathrm{ind}(\omega)$, which is defined recursively as:

1. if $\omega$ is the root, let $\mathrm{ind}(\omega) = 1$;
2. if $\omega$ is the left child of the node $\nu$, let $\mathrm{ind}(\omega) = 2 \times \mathrm{ind}(\nu)$;
3. otherwise, if $\omega$ is the right child of the node $\nu$, let $\mathrm{ind}(\omega) = 2 \times \mathrm{ind}(\nu) + 1$.

For a tree $t$, the set of level-order indices of the nodes is denoted by $\mathrm{IND}(t)$. The set $\mathrm{IND}(t)$ completely characterizes the topological structure of $t$, and will be a very useful device for proving theorems about this structure.

An important relationship between trees is the notion of a *subtree*, which is an analog of the idea of subset. A tree $s$ is called a *topological subtree* of a tree $t$ when every node in $s$ is also in $t$, that is, $\mathrm{IND}(s) \subseteq \mathrm{IND}(t)$. Moreover, if for every node $k \in \mathrm{IND}(s)$ the two trees also have the same nodal attributes, then $s$ is called an *attribute subtree* of $t$.

Also useful will be a set operations, such as union and intersection, on the topological binary tree space (i.e., when only structure is considered). For two binary trees $t_1$ and $t_2$, the tree $t$ is the *union (intersection) tree* if $\mathrm{IND}(t) = \mathrm{IND}(t_1) \cup \mathrm{IND}(t_2)$ [$\mathrm{IND}(t) = \mathrm{IND}(t_1) \cap \mathrm{IND}(t_2)$, resp.]. A horizon for our statistical analysis is provided by the union of all trees in the sample, which is called the *support tree*. This allows simplification of our analysis, because we only need to consider topological subtrees of the support tree.

The set of all topological subtrees of a given tree $t$ is called a *subtree class* and is denoted $\mathcal{T}_t$. The terminology "class" is used because each $\mathcal{T}_t$ is closed under union and intersection.

As noted in Section 1, the first major goal of statistical analysis of samples of tree-structured objects is careful definition of a centerpoint of the data set. For classical multivariate data, there are many notions of centerpoint, and even the simple concept of sample mean can be characterized in many ways. After careful extensive investigation, we have found that approaches related to the *Fréchet Mean* seem most natural. This characterizes the centerpoint as the binary tree which is the closest to all other trees in some sense (sum of squared Euclidean distances gives the sample mean in multivariate analysis). This requires a metric on the space of binary trees. Thus, the second fundamental issue is the definition of a distance between two trees. This will be developed first for the case of topology only, that is, without nodal attributes, in Section 3.2. In Section 3.3 this metric will be extended to properly incorporate attributes.

3.2. *Metric on the binary tree space without nodal attributes.* Given a tree $t$, its topological structure is represented by its set of level-order indices $\mathrm{IND}(t)$. Two trees have similar (different) topologies when their level-order index sets are similar (different, resp.). Hence, the noncommon level-order indices give an indication of the differences between two trees. Thus, for any two topological binary trees $s$ and $t$, define the metric

$$(3.1) \qquad d_I(s,t) = \sum_{k=1}^{\infty} 1\{k \in \mathrm{IND}(s) \triangle \mathrm{IND}(t)\},$$

where $\triangle$ is used to denote the symmetric set difference ($A \triangle B = (A \cap \overline{B}) \cup (\overline{A} \cap B)$, where $\overline{A}$ is the complement of $A$). Note that $d_I(s,t)$ counts the total number of nodes which show up only in either $s$ or $t$, but not both of them. Another useful view is that this metric is the smallest number of addition and deletion of nodes required to change the tree $s$ into $t$. Since $d_I$ is always an integer, it is called the *integer tree metric*, hence the subscript $I$. This will be extended to trees with attributes in Section 3.3 by adding a fractional part to this metric.

This metric can also be viewed in another way. Each binary tree can be represented as a binary string using 1 for an existent node and 0 otherwise. Since the metric $d_I$ counts differences between strings of 0s and 1s, it is just the Hamming distance from coding theory.

3.3. *Metric on the binary tree space with nodal attributes.* The integer tree metric $d_I$ captures topological structure of the tree population. In many important cases, including image analysis, the nodes of the trees contain useful attributes (numerical values, see Section 1) which also characterize important features of data objects.

The attributes contained in the node with level-order index $k$ on the tree $t$ are denoted by $(x_{tk}, y_{tk})$, where for simplicity only the case of two attributes per node is treated explicitly here. For each node, indexed by $k$, the sample mean attribute vector, $\sum_{t \in \mathcal{S}} (x_{tk}, y_{tk}) / \sum_{t \in \mathcal{S}} 1\{k \in \mathrm{IND}(t)\}$, can be assumed to be zero in the theoretical development, by subtracting the sample mean from the corresponding attribute vector of every tree which has the node $k$. Moreover, the upper bound of the absolute values of the attributes, $|x_{tk}|$ and $|y_{tk}|$, can be chosen as $\frac{\sqrt{2}}{4}$. Given any sample $\mathcal{S}$, this assumption can always be satisfied by multiplying each attribute by the scale factors $\frac{\sqrt{2}}{4}(\max_{t \in \mathcal{S}} |x_{tk}| 1\{k \in \mathrm{IND}(t)\})^{-1}$ and $\frac{\sqrt{2}}{4}(\max_{t \in \mathcal{S}} |y_{tk}| 1\{k \in \mathrm{IND}(t)\})^{-1}$. This can induce some bias in our statistical analysis, which can be partly controlled through careful choice of weights as discussed below, or by appropriate transformation of the attribute values. But this assumption is important to control the magnitude of the attribute component of the metric with respect to the topological component. The bound $\frac{\sqrt{2}}{4}$ is used because the Euclidean distance between two-dimensional vectors whose entries satisfy this bound is at most 1. For the general nodal attribute vector (e.g., the nodal attribute vectors of the blood vessel trees), a different bound will be chosen to make the attribute difference (between two trees) less than 1.

For any trees $s$ and $t$ with nodal attributes, define the new metric (Theorem 3.1 establishes that this is indeed a metric)

$$(3.2) \qquad \delta(s,t) = d_I(s,t) + f_\delta(s,t),$$

where

$$
(3.3) \quad
\begin{aligned}
f_\delta(s,t) = \Bigg[ &\sum_{k=1}^{\infty} \alpha_k((x_{sk} - x_{tk})^2 + (y_{sk} - y_{tk})^2) \\
&\times 1\{k \in \mathrm{IND}(s) \cap \mathrm{IND}(t)\} \\
&+ \sum_{k=1}^{\infty} \alpha_k(x_{sk}^2 + y_{sk}^2) 1\{k \in \mathrm{IND}(s) \backslash \mathrm{IND}(t)\} \\
&+ \sum_{k=1}^{\infty} \alpha_k(x_{tk}^2 + y_{tk}^2) 1\{k \in \mathrm{IND}(t) \backslash \mathrm{IND}(s)\} \Bigg]^{1/2}
\end{aligned}
$$

and where $\{\alpha_k\}_{k=1}^{\infty}$ is a nonnegative weight series with $\sum_{k=1}^{\infty} \alpha_k = 1$. These weights are included to allow user intervention on the importance of various nodes in the analysis (e.g., in some cases it is desirable for the root node to dominate the analysis, in which case $\alpha_1$ is taken to be relatively large). When there is no obvious choice of weights, equal weighting, $\alpha_k = \frac{1}{\#(\text{nodes appearing in the sample})}$ for nodes $k$ that appear in the sample, and $\alpha_k = 0$ otherwise, may be appropriate. All the theorems in this paper are developed for general weight sequences. But, in Section 3.6 of Wang and Marron [20] we consider some toy examples based on the *exponential weight* sequence, which gives the same weight to nodes within a level and uses an exponentially decreasing sequence across levels. In particular, the weight

$$(3.4) \qquad \alpha_k = \{2^{-(2i+1)}\}, \qquad \text{where } i = \lfloor \log_2 k \rfloor,$$

(where $\lfloor \cdot \rfloor$ denotes the greatest integer function) is used for each node on the $i$th level, $i = 0, 1, 2, \ldots$. In the analysis of the blood vessel data, different normalization of the attributes is required, because there are as many as six attributes per node. The data analyzed in Section 2 was first recentered to have 0 mean and rescaled so that the absolute value of the attributes was bounded by $\frac{1}{2\sqrt{7}}$. To more closely correspond to the original data, all of the displays in Section 2 are shown on the original scale.

The last two summations in equation (3.3) are included to avoid loss of information from those nodal attributes that are in one tree and not the other. This formulation, plus our assumption on the attributes, ensures that the second term in equation (3.2), $f_\delta$ (where "$f$" means fractional part of the metric), is at most 1.

Also, note that $f_\delta$ is a square root of a weighted sum of squares. When trees $s$ and $t$ have the same tree structure, $f_\delta(s, t)$ can be viewed as a weighted Euclidean distance. In particular, the nodal attributes of a tree $t$ can be combined into a single long vector called the *attribute vector*, denoted $\vec{v}$, for conventional statistical analysis. For an attribute subtree of $t$, the collection of attributes of the nodes of this subtree are a subvector of $\vec{v}$ which is called the *attribute subvector*.

When trees $s$ and $t$ have different tree structures, it is convenient to replace the nonexistent nodal attributes with $(0, 0)$. This also allows the nodal attributes to be combined into a single long vector, $\vec{v}$. Then $f_\delta(s, t)$ is a weighted Euclidean metric on these vectors.

For another view of $f_\delta$, rescale the entries of the vector by the square root of the weights $\alpha_k$. Then $f_\delta$ is the ordinary Euclidean metric on these rescaled vectors.

Next, Theorem 3.1 shows that $\delta$ is a metric. This requires the following assumption.

ASSUMPTION 1.    The weight $\alpha_k$ is positive and $\sum \alpha_k = 1$.

THEOREM 3.1.    *Under Assumption 1, $\delta$ is a metric on the tree space with nodal attributes.*

A sketch of the proof of Theorem 3.1 is in Section 4. The full proof is in the dissertation Wang [19], the proof of Theorem 3.1.2 in Section 3.1.

REMARK.    To understand why the Assumption 1 is critical to Theorem 3.1, consider the integer part $d_I$. While $d_I$ is a metric on topological tree space, it is not a metric on the binary tree space with nodal attributes. In particular, for any two binary trees $s$ and $t$ with the same topological structure, $d_I(s, t)$ is always equal to zero regardless of their attribute difference. Thus, $d_I$ is only a pseudo-metric on the tree space with nodal attributes. The Assumption 1 ensures that $\delta$ is a metric, not just a pseudo-metric.

3.4. *Central tree.*  In the Euclidean space $\mathbb{R}^1$, for a given data set of size $n$, there are two often-used measurements of the centerpoint, the sample mean and the sample median. Nonuniqueness for the median arises when $n$ is an even number. In this section the concepts of the sample median and the sample mean will be extended to the binary tree spaces, both with and without nodal attributes.

First, the case with no nodal attributes, that is, only topological structure, is considered. A sensible notion of centerpoint is the *median tree*, which is defined as the minimizing tree, $\arg\min_t \sum_{i=1}^n d_I(t, t_i)$, taken over all trees $t$.

This is a modification of the Fréchet mean, $\arg\min_t \sum_{i=1}^n d_I(t, t_i)^2$, which is used because it allows straightforward fast computation. This can be done using the characterization of the minimizing tree that is given in Theorem 3.2.

THEOREM 3.2.   *If a tree $s$ is a minimizing tree according to the metric $d_I$, then all the nodes of tree $s$ must appear at least $\frac{n}{2}$ times in the binary tree sample $\mathcal{S}$. Moreover, the minimizing tree $s$ (according to $d_I$) must contain all the nodes, which appear more than $\frac{n}{2}$ times, and may contain any subset of nodes that appear exactly $\frac{n}{2}$ times.*

The proof is given in Section 4.

Nonuniqueness may arise when the sample size is an even number. The *minimal median tree*, which has the fewest nodes among all the median trees, is recommended as a device for breaking any ties.

Banks and Constantine [1] independently developed essentially the same notion of central tree and this characterization of the minimizing tree, which is called the *majority rule*. We use this same terminology.

Next the case of nodal attributes is considered. Our proposed notion of centerpoint in this case is called the *median-mean tree*. It has properties similar to the sample median with respect to $d_I$ and similar to the sample mean with respect to $f_\delta$. Its tree structure complies with the majority rule and its nodal attributes can be calculated as the sample mean $\sum_{t\in\mathcal{S}}(x_{tk}, y_{tk})/\sum_{t\in\mathcal{S}} 1\{k \in \text{IND}(t)\}$. As for the median tree, the median-mean tree may not be unique, and again the minimal median-mean tree (with minimal number of nodes) is suggested for breaking such ties.

The median-mean tree is not always the same as the Fréchet mean,

$$\arg\min_t \sum_{s\in\mathcal{S}} \delta(t, s)^2.$$

We recommend the median-mean tree because it is much faster to compute. The median-mean tree is also most natural as the centerpoint of the Pythagorean theorem (i.e., sums of squares analysis) developed in Section 3.5.

Another useful concept is the *average support tree*, which consists of all the nodes that appear in the tree sample with nodal attributes calculated as averages, as done in the median-mean tree. Thus the median-mean tree is an attribute subtree of the average support tree.

3.5. *Variation analysis in the binary tree space with nodal attributes.* Now that the central tree has been developed, the next question is how to quantify the variation of the sample about the centerpoint, that is, about the median-mean tree.

In Euclidean space, the classical analysis of variance approach based on decomposing sums of squares provides a particularly appealing approach to quantifying variation. This analysis has an elegant geometric representation via the Pythagorean theorem.

After a number of trials, we found that the most natural and computable analog of the classical ANOVA decomposition came from generalizing the usual squared Euclidean norm to the *variation function*,

$$(3.5) \qquad V_\delta(s,t) = d_I(s,t) + f_\delta^2(s,t).$$

Note that if every tree has the same structure, then this reduces to classical sums of squares, and the median-mean tree is the Fréchet mean, with respect to the variation $V_\delta(s,t) = f_\delta^2(s,t)$, in the sense that it is the minimizing tree over $t$ of $\sum_{s \in \mathcal{S}} V_\delta(s,t)$. This is also true in the case of tree samples that are purely topological, that is, that have no attributes, when the variation becomes $V_\delta(s,t) = d_I(s,t)$. Then $d_I$ is a full metric (not just a pseudo-metric), which can be written as a sum of zeros and ones [see equation (3.1)]. So the metric $d_I$ can be interpreted as a sum of squares, because

$$(3.6) \quad \sum_{k=1}^{\infty} (1\{k \in \mathrm{IND}(s) \triangle \mathrm{IND}(t)\})^2 = \sum_{k=1}^{\infty} 1\{k \in \mathrm{IND}(s) \triangle \mathrm{IND}(t)\} = d_I(s,t).$$

In Euclidean space, the total variation of a sample can be measured by the sum of squared distances to its sample mean. For the tree sample $\mathcal{S}$ and the median-mean tree $m_\delta$, the total variation about the median-mean is defined as

$$\sum_{s \in \mathcal{S}} V_\delta(s, m_\delta) = \sum_{s \in \mathcal{S}} d_I(s, m_\delta) + \sum_{s \in \mathcal{S}} f_\delta^2(s, m_\delta).$$

This total variation about the median-mean tree does not depend on how the tie is broken between the median-mean trees (when it is not unique).

In classical statistics, PCA is a useful tool to capture the features of a data set by decomposing the total variation about the centerpoint. In PCA the first principal component eigenvector indicates the direction in which the data vary the most. Furthermore, other eigenvectors maximize variation in successive orthogonal residual spaces.

In binary tree space, each tree in the sample is considered to be a data point. Unlike Euclidean space, binary tree space is a nonlinear space according to the metric $\delta$ defined at (3.2). As noted above, because the space is nonlinear, the generalization of PCA is not straightforward. The foundation of our analog of PCA, is a notion of a one-dimensional manifold in binary tree space, which is a set of trees that plays the role of a "line" (a one-dimensional subspace in Euclidean space). There are two important types, defined below, both of which are called *treeline*.

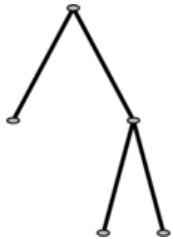

Fig. 8.  *Toy example tree $t$, for illustrating the concept of s-treeline, shown in Figure 9.*

DEFINITION 3.1.  Suppose $l = \{u_0, u_1, u_2, \ldots, u_m\}$ is a set of trees with (or without) nodal attributes in the subtree class $\mathcal{T}_t$ of a given tree $t$. The set $l$ is called a *structure treeline* (*s-treeline*) starting from $u_0$ if for $i = 1, 2, \ldots, m$,

1. $u_i$ can be obtained by adding a single node (denoted by $\nu_i$) to the tree $u_{i-1}$ (thus, when attributes exist, they are common through the treeline);
2. the next node to be added, $\nu_{i+1}$, is the child of $\nu_i$;
3. the first tree $u_0$ is *minimal*, in the sense that the ancestor node of $\nu_1$ is the root node, or else has another child.

REMARK.  Structure treelines are "one-dimensional" in the sense that they follow a single path, determined by $u_0$ and the sequence of added nodes $\nu_1, \ldots, \nu_m$. In this sense the elements of $l$ are *nested*. Also when there are attributes, each attribute vector is the corresponding attribute subvector (defined in Section 3.3) of its successor.

In Definition 3.1 the tree $u_{i-1}$ is a subtree (an attribute subtree, if there are attributes) of the trees $u_i$, $u_{i+1}$, and so on. Since every element in the $s$-treeline is a topological subtree of $t$, the length of the $s$-treeline cannot exceed the number of levels of the tree $t$. Illustration of the concept of $s$-treeline is shown in Figures 8 and 9.

Figure 9 shows the $s$-treeline in $\mathcal{T}_t$, where $t$ has the tree structure shown in Figure 8. Figure 9 indicates both tree topology and attributes. The positive attributes $(x, y)$ are graphically illustrated with a box for each node, where $x$ is shown as the horizontal length and $y$ is the height.

Note that the attributes are common for each member of the treeline. Each succeeding member comes from adding a new node. The starting member, $u_0$, cannot be reduced, because the root node has a child which does not follow the needed sequence.

A structure treeline $l$ is said to *pass through* the tree $u$, when the tree $u$ is an element of the tree set $l$, that is, $u \in l$. Recall from Section 2 that, for the blood vessel data, Figure 5 shows the topology of the structure treeline passing through the median-mean tree, and Figure 6 shows the corresponding attributes. The central tree in each figure is the median-mean tree.

An $s$-treeline summarizes a direction of changing tree structures. The following definition will describe a quite different direction in tree space, in which all trees have the same tree structure but changing nodal attributes.

DEFINITION 3.2. Suppose $l = \{u_\lambda : \lambda \in \mathbb{R}\}$ is a set of trees with nodal attributes in the subtree class $\mathcal{T}_t$ of a given tree $t$. The set $l$ is called an *attribute treeline* ($a$-treeline) *passing through* a tree $u^*$ if

1. every tree $u_\lambda$ has the same tree structure as $u^*$;
2. the nodal attribute vector is equal to $\vec{v}^* + \lambda \vec{v}$, where $\vec{v}^*$ is the attribute vector of the tree $u^*$ and where $\vec{v}$ is some fixed vector, $\vec{v} \neq \vec{0}$.

REMARK. An $a$-treeline is determined by the tree $u^*$ and the vector $\vec{v}$. The treeline is "one-dimensional" in this sense, which is essentially the same as a line in Euclidean space.

Figure 10 shows some members of an $a$-treeline from the same subtree class $\mathcal{T}_t$ shown in Figure 8, with $\lambda = 0.5, 1.0, 1.2, 1.5$ and

$$\vec{v} = [0.2, 0.1, 0.1, 0.2, 0.1, 0.1, 0.2, 0.2]'.$$

The topological structures of all of the trees in Figure 10 are the same. The dimensions of the boxes, illustrating the values of the attributes, change linearly.

In Section 2, Figure 7 illustrated an attribute treeline. That treeline highlighted the strong variation between orientations of the trees in the sample.

From now on, both $s$-treelines and $a$-treelines are called treelines. An analogy of the first principal component is the treeline which explains most of the variation in the data. A notion of *projection*, of a tree onto a treeline, needs to be defined, because this provides the basis for decomposition of sums of squares.

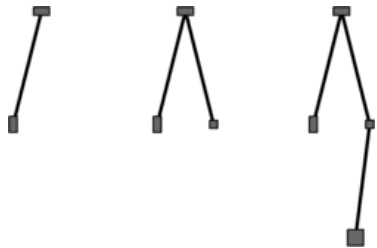

FIG. 9. *An s-treeline in $\mathcal{T}_t$, based on the tree $t$ shown in Figure 8.*

For any tree $t$ and treeline $l$, the projection of $t$ onto $l$, denoted $P_l(t)$, is the tree which minimizes the distance $\delta(t, \cdot)$ over all trees on the treeline $l$. The idea of projection is most useful when it is unique, as shown in the next theorem.

THEOREM 3.3. *Under Assumption 1, the projection of a tree $t$ onto a treeline $l$ is unique.*

Because of an editorial decision to save space, the proof is given in Section 4 of Wang and Marron [20].

The Pythagorean theorem is critical to the decomposition of the sums of squares in classical ANOVA. Analogs of this are now developed for tree samples. Theorem 3.4 gives a Pythagorean theorem for $a$-treelines and Theorem 3.5 gives a Pythagorean theorem for $s$-treelines.

THEOREM 3.4 (Tree version of the Pythagorean theorem: Part I). *Let $l$ be an $a$-treeline passing through a tree $u$ in the subtree class $\mathcal{T}_t$. Then, for any $t \in \mathcal{T}_t$,*

$$(3.7) \qquad V_\delta(t, u) = V_\delta(t, P_l(t)) + V_\delta(P_l(t), u).$$

REMARK. This states that the variation (our analog of squared distance) of a given tree $t$ from a tree $u$ in the treeline $l$, which is essentially the hypotenuse of our triangle, is the sum of the variation of $t$ from $P_l(t)$, plus the variation of $P_l(t)$ from $u$, representing the legs of our triangle. This is the key to finding treelines that explain maximal variation in the data, because $V_\delta(t, u)$ is independent of $l$, so maximizing (over treelines $l$) a sample sum over $V_\delta(P_l(t), u)$ is equivalent to minimizing the residual sum over $V_\delta(t, P_l(t))$.

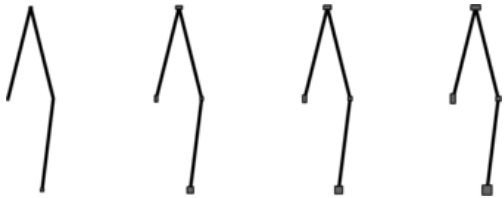

FIG. 10. *Toy example of an $a$-treeline for the same subtree class $\mathcal{T}_t$ as in Figure 8. Several members of the treeline are shown. The attributes are a linear function of each other.*

In this paper, only those $s$-treelines where every element is an attribute subtree of the average support tree (as defined in Section 3.4) are considered, because this gives a tree version of the Pythagorean theorem, shown next.

THEOREM 3.5 (Tree version of the Pythagorean theorem: Part II).  *Let* $\mathcal{S} = \{t_1, t_2, \ldots, t_n\}$ *be a sample of trees. Let* $l$ *be an* $s$-*treeline where every element is an attribute subtree of the average support tree of* $\mathcal{S}$. *Then, for any* $u \in l$,

$$(3.8) \qquad \sum_{t \in \mathcal{S}} V_\delta(t, u) = \sum_{t \in \mathcal{S}} V_\delta(t, P_l(t)) + \sum_{t \in \mathcal{S}} V_\delta(P_l(t), u).$$

REMARK.  This theorem complements Theorem 3.4, because it now gives a structure treeline version of the Pythagorean theorem, which simplifies analysis of variance, because minimizing the residual sum $\sum_{t \in \mathcal{S}} V_\delta(P_l(t), t)$ is equivalent to maximizing the sum $\sum_{t \in \mathcal{S}} V_\delta(\mu_\delta, P_l(t))$ over all treelines passing through the minimal median-mean tree, $\mu_\delta$. In some sense, this theorem is not so strong as Theorem 3.4, because the sample summation is needed, while the Pythagorean theorem 3.4 is true even term by term.

Sketches of the proofs of Theorems 3.4 and 3.5 are given in Section 4 of Wang and Marron [20]. Further details are in the proofs of Theorems 3.5.3 and 3.5.4 in Section 3.5 of Wang [19].

The foundations are now in place to develop variation analysis in binary tree space. There are two main steps to the *PCA on trees* variation analysis.

First, find an $s$-treeline $l_{\mathrm{PS}}$ that minimizes the sum $\sum_{t \in \mathcal{S}} V_\delta(t, P_l(t))$ over $l$ passing through the minimal median-mean tree $\mu_\delta$ of the sample $\mathcal{S}$, that is,

$$(3.9) \qquad l_{\mathrm{PS}} = \arg\min_{l \,:\, \mu_\delta \in l} \sum_{t \in \mathcal{S}} V_\delta(t, P_l(t)).$$

This structure treeline is called a *one-dimensional principal structure representation* (*treeline*) of the sample $\mathcal{S}$. Because of the Pythagorean theorem 3.5, the one-dimensional structure treeline $l_{\mathrm{PS}}$ explains a maximal amount of the variation in the data, as is done by the first principal component in Euclidean space. This is illustrated in the context of the blood vessel data in Section 2. Figure 6 shows the principal structure treeline $l_{\mathrm{PS}} = \{u_0, u_1, u_2\}$ with nodal attributes, where $u_1$ is the unique median-mean tree (also the minimal median-mean tree) of the sample. Figure 5 shows the topological tree structures of the principal structure treeline in Figure 6.

Second, a notion of principal attribute treeline direction will be developed. This will complement the principal structure treeline, in the sense that together they determine an analog of a two-dimensional subspace of binary tree space. Recall from Definition 3.2 that an attribute treeline is indexed by a starting tree $u^*$, with attribute vector $\vec{v}^*$, and by a direction vector $\vec{v}$, and has general attribute vector $\vec{v}^* + \lambda \vec{v}$, for $\lambda \in \mathbb{R}$. To create the desired two-dimensional structure, we consider a family of attribute treelines indexed by the nested members $\{u_0, u_1, \ldots, u_m\}$ of the principal structure treeline and their corresponding nested (in the sense of attribute subvectors, as defined in Section 3.3) attribute vectors $\{\vec{v}_0^*, \vec{v}_1^*, \ldots, \vec{v}_m^*\}$, and indexed by a set of nested direction vectors $\{\vec{v}_0, \vec{v}_1, \ldots, \vec{v}_m\}$.

The union of treelines that are nested in this way is called a *family of attribute treelines*. This concept is developed in general in the following definition.

DEFINITION 3.3. Let $l = \{u_0, u_1, \ldots, u_m\}$ be a structure treeline, and let $\vec{c}$ be a vector of attributes corresponding to the nodes of $u_m$. The $l, \vec{c}$-*induced family of attribute treelines*, $\mathcal{E}_{l,\vec{c}} = \{e_0, e_1, \ldots, e_m\}$, is defined for $k = 0, 1, \ldots, m$ as

$$e_k = \{t_\lambda : t_\lambda \text{ has attribute vector } \vec{v}_k + \lambda \vec{c}_k, \ \lambda \in \mathbb{R}\},$$

where $\vec{v}_k$ is the attribute vector of $u_k$, and where $\vec{c}_k$ is the corresponding attribute subvector of $\vec{c}$.

Next an appropriate family of attribute treelines is chosen to provide maximal approximation of the data (as is done by the first two principal components in Euclidean space). Following conventional PCA, we start with the principal structure treeline $l_{\mathrm{PS}}$ (which we will denote in this paragraph as $l$ simply to save a level of subscripting) and choose the direction vector $\vec{c}$ so that the $l, \vec{c}$-induced family of attribute treelines explains as much of the data as possible. In particular, we define the principal attribute direction vector, $\vec{c}_{\mathrm{PA}}$, as

$$\vec{c}_{\mathrm{PA}} = \arg\min_{\vec{c}} \sum_{t \in \mathcal{S}} V_\delta(t, P_{\mathcal{E}_{l,\vec{c}}}(t)),$$

where $P_{\mathcal{E}_{l,\vec{c}}}(t)$ is the projection of the tree $t$ onto the $l, \vec{c}$-induced family of attribute treelines. This is an analog of a two-dimensional projection, defined as

$$P_{\mathcal{E}_{l,\vec{c}}}(t) = \arg\min_{s \, : \, s \in e(t)} \delta(t, s),$$

where $e(t)$ is the attribute treeline determined by the tree $P_l(t)$ (the projection of $t$ onto the principal structure treeline) and by the direction vector which is the corresponding attribute subvector of $\vec{c}$.

The elements of the $l_{\mathrm{PS}}, \vec{c}_{\mathrm{PA}}$-induced family of attribute treelines are all called *principal attribute treelines*. As the first two principal components can illuminate important aspects of the variation in data sets of curves, as demonstrated by, for example, Ramsay and Silverman [16, 17], the principal structure and attribute treelines can find important structure in a data sets of trees. An example of this is shown in Figure 7 of this paper and Figure 10 of Wang and Marron [20], where the principal attribute treeline through the median-mean tree revealed the change in orientation in the data.

For completely different extensions of PCA in nonlinear ways, see the principal curve idea of Hastie and Stuetzle [9] and the principal geodesic analysis of Fletcher, Lu and Joshi [8]. Principal curves provide an interesting nonlinear decomposition of Euclidean space. Principal geodesics provide a decomposition of data in nonlinear manifolds, including Lie groups and symmetric spaces.

For an interesting comparison of our tree version PCA and regular PCA, see Section 3.6 of Wang and Marron [20].

## 4. Derivations of theorems.

A SKETCH OF THE PROOF OF THEOREM 3.1.   The proof follows from the fact that $d_I$ is a metric on the binary tree space without nodal attributes and $f_\delta$ is a weighted Euclidean distance between two attribute vectors (the nodal attributes for nonexistent nodes are treated as zeros).   □

PROOF OF THEOREM 3.2.   Let $s$ be a minimizing tree according to the integer tree metric $d_I$. Suppose some of the nodes in $s$ appear less than $\frac{n}{2}$ times and $\nu$ is the node with the largest level among all of those nodes. If a node appears less than $\frac{n}{2}$ times, so do its children. Thus, $\nu$ must be a terminal node of $s$.

For the binary tree $s'$, with $\mathrm{IND}(s') = \mathrm{IND}(s) \backslash \{\mathrm{ind}(\nu)\}$, the following equation is satisfied:

$$(4.1) \qquad \sum_{i=1}^n d_I(s', t_i) = \sum_{i=1}^n d_I(s, t_i) + n_\nu - (n - n_\nu),$$

where $n_\nu = \#\{$appearance of the node $\nu$ in the sample $\mathcal{S}\}$. Since $n_\nu < \frac{n}{2}$,

$$\sum_{i=1}^n d_I(s', t_i) < \sum_{i=1}^n d_I(s, t_i),$$

which is a contradiction of the assumption that $s$ is a minimizing tree.

From the proof above, if $n_\nu = \frac{n}{2}$, then $\sum_{i=1}^n d_I(s', t_i) = \sum_{i=1}^n d_I(s, t_i)$; that is, $s'$ is also a minimizing tree. Therefore, the minimizing tree may contain any subset of the nodes that appear exactly $\frac{n}{2}$ times.

Finally, a proof is given of the fact that the minimizing binary tree $s$ contains all the nodes which appear more than $\frac{n}{2}$ times.

Suppose the node $\omega$ appears more than $\frac{n}{2}$ times in the sample $\mathcal{S}$ and $\mathrm{ind}(\omega) \notin \mathrm{IND}(s)$. Without loss of generality, suppose that $\omega$ is a child of some node in the binary tree $s$. Otherwise, choose one of its ancestor nodes.

For the binary tree $s''$, with $\text{IND}(s'') = \text{IND}(s) \cup \{\text{ind}(\omega)\}$, the following equation is satisfied:

$$(4.2) \qquad \sum_{i=1}^{n} d_I(s, t_i) = \sum_{i=1}^{n} d_I(s'', t_i) + n_\omega - (n - n_\omega),$$

where $n_\omega = \#\{\text{appearance of the node } \omega \text{ in the sample } \mathcal{S}\}$. Since $n_\omega > \frac{n}{2}$,

$$\sum_{i=1}^{n} d_I(s'', t_i) < \sum_{i=1}^{n} d_I(s, t_i),$$

which is a contradiction of the assumption that $s$ is the minimizing tree. $\square$

**Acknowledgments.** Much of this paper is the dissertation work of the first author, written under the supervision of the second. S. M. Pizer and the UNC MIDAG project provided financial support for the both authors, and many helpful conversations throughout this research. The term "object oriented data analysis" came from a conversation with J. O. Ramsay. E. Bullitt provided the blood vessel tree data (gathered with support from grant R01 EB000219 NIH-NIBIB), and the viewer STree that was used to construct Figure 2. Paul Yushkevich provided the program mmview that made Figure 1.

Department of Statistics
Colorado State University
Fort Collins, Colorado 80523
USA
E-mail: wanghn@stat.colostate.edu

Department of Statistics
University of North Carolina
Chapel Hill, North Carolina 27599
USA
E-mail: marron@email.unc.edu